\documentclass[12pt]{amsart}
\usepackage{amsmath,amsthm} 
\usepackage{amssymb,amsfonts,amsxtra}
\usepackage{psfrag,graphicx}
\usepackage{hyperref}

\theoremstyle{plain}
\newtheorem{theorem}{Theorem}
\newtheorem{proposition}[theorem]{Proposition}
\newtheorem{corollary}[theorem]{Corollary}
\newtheorem{lemma}[theorem]{Lemma}

\theoremstyle{definition}
\newtheorem{definition}[theorem]{Definition}
\newtheorem{remark}[theorem]{Remark}

\newtheorem*{claim}{Claim}

\numberwithin{theorem}{section}

\makeatletter
  \def\tagform@#1{\maketag@@@{%
   \textbf{(\ignorespaces#1\unskip\@@italiccorr)}}}%
   \renewcommand{\eqref}[1]{\textup{\maketag@@@{(\ignorespaces%
        {\ref{#1}}\unskip\@@italiccorr)}}}
\let\@itemize@\itemize
\def\itemize{\@itemize@\parskip 0pt\relax}
\def\@listi{\leftmargin28.5pt\parsep 0pt\topsep 2pt plus1pt minus1pt
 \itemsep2pt plus2pt minus1pt}
\let\@listI\@listi
\@listi
\makeatother

\newcommand{\zindex}[3]{\put(#1,#2){\makebox(0,0){${#3}$}}}
\newenvironment{pict}[2]%
	{\setlength{\unitlength}{1mm}
	\begin{center}
	\begin{picture}(#1,#2)
	\scriptsize
}%
	{\end{picture}
	\end{center}
 	\noindent}

\DeclareFontFamily{OML}{rsfs}{\skewchar\font'177}
\DeclareFontShape{OML}{rsfs}{m}{n}{ <5> <6> rsfs5 <7> <8> <9> rsfs7
  <10> <10.95> <12> <14.4> <17.28> <20.74> <24.88> rsfs10 }{}
\DeclareMathAlphabet{\mathfs}{OML}{rsfs}{m}{n}

\newcommand{\As}{{\mathfs A}}

\newcommand{\co}{\colon\thinspace}

\newcommand{\DS}{\ensuremath{{\mathfs D}}} 

\renewcommand{\epsilon}{\varepsilon}

\newcommand{\iv}{o} 
\newcommand{\tv}{t}

\DeclareMathOperator{\Aut}{Aut}

\begin{document}


\title{Whitehead moves for $G$--trees}
\author{Matt Clay and Max Forester}
\address{Mathematics Department, 
        University of Oklahoma, 
        Norman, OK 73019, 
        USA}
\email{mclay@math.ou.edu\\forester@math.ou.edu}

\begin{abstract}
We generalize the familiar notion of a Whitehead move from Culler
and Vogtmann's Outer space to the setting of deformation spaces of
$G$--trees.  Specifically, we show that there are two moves, each of
which transforms a reduced $G$--tree into another reduced $G$--tree,
that suffice to relate any two reduced trees in the same
deformation space.  These two moves further factor into three moves
between reduced trees that have simple descriptions in terms of graph of
groups.  This result has several applications. 
\end{abstract}

\maketitle

\thispagestyle{empty}


\section{Introduction}\label{sc:intro}

Whitehead automorphisms of $F_{n}$ (the free group of rank $n$) 
generate the automorphism group $\Aut(F_{n})$.  These automorphisms were
used by J.H.C. Whitehead to construct an algorithm to decide whether
two given elements of $F_{n}$ are related by an automorphism 
\cite{ar:W36}. 

These automorphisms can be interpreted as certain moves on an
$n$--rose whose fundamental group is marked with an isomorphism to
$F_{n}$ \cite{ar:H79}.  As such, they can be used to provide a path of
marked $n$--roses connecting any two marked $n$--roses in
Culler and Vogtmann's Outer space \cite{ar:CV86}.  This is
the space of marked metric graphs modulo homothety. 
By passing to the universal covers of the marked metric
graphs, an alternative description of Outer space is the space of free
minimal actions of $F_{n}$ on metric simplicial trees, again modulo
homothety. 

Deformation spaces of $G$--trees (see \cite{ar:F02}) are a generalization
of Outer space, where the actions of a group $G$ on a simplicial tree are
allowed to have nontrivial stabilizers, but the set of elliptic subgroups
(subgroups that fix
points) is uniform throughout the deformation space.  See
\cite{ar:Cpp,ar:C05,ar:F06,ar:GL07,ar:GL07a,ar:L07,ar:MM} for examples and
applications of deformation spaces.

In a deformation space, the analogue of a rose in Outer space is a
\emph{reduced} tree (defined in Section \ref{sc:def}). 
The purpose of this note is to find a finite set of moves, analogous
to Whitehead moves, that will provide a path of reduced $G$--trees
through a deformation space, connecting any two given reduced $G$--trees.
This is achieved in Theorem \ref{th:whm}, where it is shown that two
particular moves suffice. These two
moves are then decomposed into three simpler moves called \emph{slide},
\emph{induction}, and \emph{$\As^{\pm 1}$--moves}. Our main theorem is
the following. 

\begin{theorem}\label{th:moves}
In a deformation space of $G$--trees, any two reduced trees are related
by a finite sequence of slides, inductions, and $\As^{\pm 1}$--moves,
with all intermediate trees reduced. 
\end{theorem}

In this paper, sliding a collection of edges over one edge is
considered a single slide move. Of course, in the case of cocompact
$G$--trees, the theorem is still valid using the traditional definition
(sliding one edge orbit at a time); see Definition \ref{def:slide}. 

One immediate consequence of Theorem \ref{th:moves} is a
strengthened form of the uniqueness theorem for JSJ decompositions of
finitely generated groups \cite{ar:F03}. Here, JSJ decompositions are
meant in the sense of Rips  and Sela \cite{ar:RS97}, Dunwoody and Sageev
\cite{ar:DSa99}, or Fujiwara and Papasoglu \cite{ar:FP06}.  

\begin{corollary}
Any two JSJ decompositions of a finitely generated group are related by a
finite sequence of slides, inductions, and $\As^{\pm 1}$--moves between
reduced decompositions. 
\end{corollary}

At the end of the paper we discuss two further applications. 
One observation is that if the deformation space is \emph{non-ascending}
(see below) then induction and $\As^{\pm 1}$--moves cannot occur. Thus
any two reduced trees are related by slide moves. This result has 
previously appeared as \cite[Theorem 7.4]{ar:F06} and
\cite[Theorem~7.2]{ar:GL07}, and indeed our proof of Theorem
\ref{th:moves} is similar in spirit to the proof given in \cite{ar:GL07}. 
The theorem also 
directly implies 
the rigidity theorem for $G$--trees \cite{ar:F02,ar:G03}, in its most
general form due to Levitt \cite{ar:L05a}. 

Lastly, 
Theorem \ref{th:moves} plays a significant role
in the solution to the isomorphism problem for certain 
generalized Baumslag--Solitar groups. This work appears in \cite{ar:CF2}.


\section{Deformation spaces}\label{sc:def}


A graph $\Gamma$ is given by $(V(\Gamma),E(\Gamma),\iv,\tv,\bar{
\ })$ where $V(\Gamma)$ are the vertices, $E(\Gamma)$ are the oriented
edges, $\iv,\tv \co E(\Gamma) \to V(\Gamma)$ are the originating and
terminal vertex maps and $\bar{ \ }\co E(\Gamma) \to E(\Gamma)$ is a
fixed point free involution, which reverses the orientations of edges.
An \emph{edge path} $\gamma = (e_{0},\ldots,e_{k})$ is a sequence of
edges such that $\tv(e_{i}) = \iv(e_{i+1})$ for $i = 0,\ldots,k-1$.  A
\emph{loop} is an edge $e \in E(\Gamma)$ such that $\iv(e) = \tv(e)$. A
\emph{geometric   edge} is a pair of the form $\{e, \bar{e}\}$. When we
say that $e, f$ are ``distinct geometric edges'' we mean that none of the
oriented edges $e, \bar{e}, f, \bar{f}$ coincide. 

Let $G$ be a group.  A \emph{$G$--tree} is a simplicial tree $T$ together
with an action of $G$ by simplicial automorphisms, without inversions 
(that is, $ge \neq \bar{e}$ for all $g \in G, e \in E(T)$).  Two
$G$--trees are considered equivalent if there is a $G$--equivariant
isomorphism between them. The quotient graph $T/G$ has the structure of a
graph of groups with a marking (an identification of $G$ with the
fundamental group of the graph of groups). We call such graphs
\emph{marked graphs of groups}, or \emph{marked graphs} for short. 

Given a $G$--tree $T$, a subgroup $H \subseteq G$ is \emph{elliptic}
if it fixes a point of $T$.  There are two moves one can perform
on a $G$--tree without changing the elliptic subgroups, called
\emph{collapse and expansion moves}; they correspond to the natural
isomorphism $A \ast_{B} B \cong A$. The exact definition is as follows.

\begin{definition}\label{def:elementarymoves}
An edge $e$ in a $G$--tree $T$ is \emph{collapsible} if $G_{e} =
G_{\iv(e)}$ and its endpoints are not in the same orbit. If one collapses
$\{e, \bar{e}\}$ and all of its translates to vertices, the resulting
$G$--tree is said to be obtained from $T$ by a \emph{collapse move}. The
reverse of this move is called an \emph{expansion move}. 

If $\Gamma$ is the corresponding marked graph of groups, $e\in
E(\Gamma)$ is collapsible if it is not a loop and the inclusion map $G_e
\hookrightarrow G_{\iv(e)}$ is surjective. The marked graph obtained from
$\Gamma$ by collapsing $e$ is denoted $\Gamma_e$. If $F \subset \Gamma$
is a forest whose edges can be collapsed iteratively, we denote the
resulting marked graph $\Gamma_{F}$. A \emph{non-trivial} forest is a
forest containing at least one edge. 

A $G$--tree (or marked graph) is \emph{reduced} if it does not admit a
collapse move. An \emph{elementary deformation} is a finite sequence of
collapse and expansion moves. Given a $G$--tree $T$, the \emph{deformation
  space} $\DS$ of $T$ is the set of all $G$--trees related to $T$ by an
elementary deformation. If $T$ is cocompact then $\DS$ is equivalently the
set of all $G$--trees having the same elliptic subgroups as $T$
\cite{ar:F02}. Equivalently, $\DS$ may be thought of as the space of
marked graphs related to $T/G$ by collapse and expansion moves. 
\end{definition}


There are three special deformations that will be considered as basic
moves. We define them below in terms of graphs of groups, but first we
need some terminology. Suppose a graph of groups has an edge $e$ which
is a loop. Let $A$ be the vertex group and $B$ the edge group, with
inclusion maps $i_0, i_1 \co B \hookrightarrow A$. If one of these maps,
say $i_0$, is an isomorphism, then $e$ is an \emph{ascending loop}. The
\emph{monodromy} is the composition $i_1 \circ i_0^{-1} \co A
\hookrightarrow A$. If the monodromy is not surjective then $e$ is a
\emph{strict ascending loop}. A deformation space $\DS$ is
\emph{ascending} if it contains a $G$--tree whose quotient graph of
groups has a strict ascending loop. Otherwise it is called
\emph{non-ascending}. 


\begin{definition}\label{def:slide}
The deformation shown below 
is called a \emph{slide move}. 
The edge groups of the edges that slide do not change. However, 
in order to perform the move, these edge groups must be contained in $C$
(considered as subgroups of $A$ before the move). 
\begin{pict}{120}{13}
\thicklines
\put(102,4){\circle*{1}}
\put(114,4){\circle*{1}}
\put(102,4){\line(1,0){12}}
\put(114,4){\line(-1,2){4}}
\put(114,4){\line(-4,3){7}}
\put(114,4){\line(0,1){8.5}}

\zindex{102}{1.5}{A}
\zindex{114}{1.5}{B}
\zindex{108}{2}{C}

\thinlines
\put(114,4){\line(5,3){5}}
\put(114,4){\line(5,-3){5}}
\put(102,4){\line(-5,3){5}}
\put(96,4){\line(1,0){6}}
\put(102,4){\line(-5,-3){5}}


\thicklines
\put(6,4){\circle*{1}}
\put(18,4){\circle*{1}}
\put(6,4){\line(1,0){12}}
\put(6,4){\line(1,2){4}}
\put(6,4){\line(0,1){8.5}}
\put(6,4){\line(4,3){7}}

\zindex{6}{1.5}{A}
\zindex{18}{1.5}{B}
\zindex{12}{2}{C}

\thinlines
\put(18,4){\line(5,3){5}}
\put(18,4){\line(5,-3){5}}
\put(6,4){\line(-5,3){5}}
\put(0,4){\line(1,0){6}}
\put(6,4){\line(-5,-3){5}}


\thicklines
\put(50,4){\circle*{1}}
\put(60,4){\circle*{1}}
\put(70,4){\circle*{1}}
\put(50,4){\line(1,0){20}}
\put(60,4){\line(0,1){8}}
\put(60,4){\line(1,2){3.6}}
\put(60,4){\line(-1,2){3.6}}

\zindex{50}{1.5}{A}
\zindex{55}{2}{C}
\zindex{60}{1.5}{C}
\zindex{65}{2}{C}
\zindex{70}{1.5}{B}

\thinlines
\put(70,4){\line(5,3){5}}
\put(70,4){\line(5,-3){5}}
\put(50,4){\line(-5,3){5}}
\put(44,4){\line(1,0){6}}
\put(50,4){\line(-5,-3){5}}


\zindex{34}{6.5}{\mbox{exp.}}
\zindex{86}{7}{\mbox{coll.}}
\put(29,5){\vector(1,0){10}}
\put(81,5){\vector(1,0){10}}
\end{pict}%

The set of edges that slide may have any cardinality. If this cardinality
is finite, however, then the edges may of course be slid one at a
time. Notice that in this situation, if the initial and final marked
graphs are reduced, then so are the intermediate graphs when edges are
slid separately. 

The edge carrying $C$ is allowed to be a loop. In this case
the only change to the graph of groups is in the inclusion maps of the
edge groups to $A$. Specifically, if $i_0, i_1\co C \hookrightarrow A$
are the inclusion maps of the loop, and $j\co D \hookrightarrow A$ is the
inclusion map of an edge, with $j(D) \subset i_0(C)$, the map $j$ is
replaced by $i_1 \circ i_0^{-1} \circ j$; see \cite[Section
3.6]{ar:F02}. Note that this results in a new marking of the graph of
groups, even if the underlying graph is unchanged. 
\end{definition}

\begin{definition}\label{def:ind} 
An \emph{induction move} is an expansion and collapse along an ascending
loop. In the diagram below the ascending loop has vertex group $A$ and
monodromy $\phi \co A \hookrightarrow A$, and $B$ is a subgroup such that
$\phi(A) \subseteq B \subseteq A$. The map $\iota\co B \hookrightarrow A$
is inclusion. The lower edge is expanded and the upper edge is collapsed,
resulting in an ascending loop with monodromy the induced map
$\phi\vert_{B}\co B \hookrightarrow B$.  

\begin{pict}{104}{12}
\thicklines
\put(92,6){\circle{10}}
\put(97,6){\circle*{1}}

\put(51,6){\oval(10,10)[b]}

\thinlines
\put(9,6){\circle{10}}
\put(14,6){\circle*{1}}

\put(51,6){\circle{10}}
\put(46,6){\circle*{1}}
\put(56,6){\circle*{1}}

\put(14,6){\line(1,1){4}}
\put(14,6){\line(1,-1){4}}

\put(56,6){\line(1,1){4}}
\put(56,6){\line(1,-1){4}}

\put(97,6){\line(1,1){4}}
\put(97,6){\line(1,-1){4}}

\scriptsize
\zindex{5.7}{6}{\phi}
\qbezier(4,3)(2.3,6)(4,9)
\put(4,3){\vector(2,-3){0}}
\zindex{16.7}{6.2}{A}

\put(25,6){\vector(1,0){12}}
\zindex{31}{7.5}{\mbox{exp. }}

\qbezier(48,11)(51,12.7)(54,11)
\put(48,11){\vector(-3,-2){0}}
\zindex{51}{9}{\phi}
\qbezier(48,1)(51,-0.7)(54,1)
\put(54,1){\vector(3,2){0}}
\zindex{51}{2.7}{\iota}
\zindex{43.5}{6}{B}
\zindex{58.7}{6.2}{A}

\put(67,6){\vector(1,0){12}}
\zindex{73}{8}{\mbox{coll. }}

\zindex{90.4}{6}{\phi\vert_B}
\qbezier(87,3)(85.3,6)(87,9)
\put(87,3){\vector(2,-3){0}}
\zindex{100}{6}{B}
\end{pict}%
The reverse of this move is also considered an induction move. Notice
that the vertex group may change, in contrast with slide moves. 
\end{definition}


\begin{definition}\label{def:A}
An \emph{$\As^{-1}$--move} is an induction followed by a collapse as
shown below. The move is always non-trivial, and it has some 
requirements: the loop is an ascending loop with monodromy $\phi$, 
with $\phi(A) \subseteq B \subseteq A$ (so that the induction can be 
performed), $B$ is a proper subgroup of both $A$ and $C$, and there are
no other edges incident to the loop. 

\begin{pict}{120}{12}
\thicklines
\put(7,6){\circle{10}}
\put(12,6){\line(1,0){13}}

\put(58,6){\circle{10}}
\put(63,6){\line(1,0){13}}

\put(109,6){\circle{10}}
\put(114,6){\circle*{1}}

\thinlines
\put(12,6){\circle*{1}}
\put(25,6){\circle*{1}}
\put(25,6){\line(1,1){4}}
\put(25,6){\line(1,-1){4}}

\put(63,6){\circle*{1}}
\put(76,6){\circle*{1}}
\put(76,6){\line(1,1){4}}
\put(76,6){\line(1,-1){4}}

\put(114,6){\line(1,1){4}}
\put(114,6){\line(1,-1){4}}

\scriptsize
\zindex{3.7}{6}{\phi}
\qbezier(2,3)(0.3,6)(2,9)
\put(2,3){\vector(2,-3){0}}
\zindex{13.5}{8.8}{A}
\zindex{18.8}{8.4}{B}
\zindex{24}{8.8}{C}

\zindex{56.3}{6}{\phi\vert_B}
\qbezier(53,3)(51.3,6)(53,9)
\put(53,3){\vector(2,-3){0}}
\zindex{64.5}{8.8}{B}
\zindex{69.8}{8.4}{B}
\zindex{75}{8.8}{C}

\zindex{117.5}{6}{C}

\put(34,6){\vector(1,0){12}}
\zindex{40}{8}{\mbox{ind.}}
\put(86,6){\vector(1,0){12}}
\zindex{92}{8}{\mbox{coll. }}
\end{pict}%

Note that before the move, the loop is a strict ascending loop,
and after, the loop is not ascending. Thus an $\As^{-1}$--move
removes an ascending loop, and its reverse, called an \emph{$\As$--move},
adds one. 

If $e$ is a loop labeled by the group $B$ with inclusion maps $i_0$ and
$i_1$, an $\As$--move can be performed on $e$ if the following
criteria are met: $i_0(B) \subsetneq i_1(B)$ or $i_1(B) \subsetneq i_0(B)$,
and both $i_0(B)$ and $i_1(B)$ are proper subgroups of the vertex group.
\end{definition}

\begin{remark}
$\As^{\pm 1}$--moves preserve the property of being reduced. The same is
not always true of slide or induction moves, unless one is in a
non-ascending deformation space. See \cite[Example 3.2]{ar:F06} for the
case of an induction move; an obvious modification also yields a slide move
example. Note as well that an $\As^{\pm 1}$--move can only
occur in an ascending deformation space. 
\end{remark}

\begin{remark}
In special cases, a move may result in a $G$--tree equivalent to the
original one. If there is a loop whose inclusion maps are isomorphisms
and are equal, then sliding all other edges incident to that vertex over
the loop will result in the same $G$--tree. For slide moves, this is the
only such example. In an induction move, the new $G$-tree is always
different, unless the quotient graph of groups is a single ascending
loop. See \cite[Theorem 2]{ar:L05a} for a thorough description of the
possibilities in this case. An $\As^{\pm 1}$--move is always non-trivial
(as are collapse and expansion moves) since the quotient graph changes. 
\end{remark}


An $n$--rose in Outer space may be regarded as a marked graph of trivial
groups with a single vertex. From this point of view, a Whitehead move is
an expansion (as defined in \ref{def:elementarymoves}) 
followed by a collapse of an edge other than the expanded one. 
The next definition generalizes this move.  Let $\DS$ be a deformation
space of a $G$--tree.  

\begin{definition}\label{def:whitehead}
Let $\Gamma,\Gamma' \in \DS$ be reduced marked graphs.  We
say $\Gamma$ and $\Gamma'$ are \emph{related by a type I Whitehead
move} if there is a marked graph $\Gamma'' \in \DS$ such that $\Gamma
= \Gamma''_{e}$ and $\Gamma' = \Gamma''_{e'}$ for distinct geometric
edges $e,e' \in E(\Gamma'')$.  We say $\Gamma$ and $\Gamma'$ are
\emph{related by a type II Whitehead move} if there is a marked graph
$\Gamma'' \in \DS$ such that $\Gamma = \Gamma''_{e}$ and $\Gamma' =
\Gamma''_{e' \cup f'}$ for distinct geometric edges $e,e',f' \in
E(\Gamma'')$. Note that Whitehead moves are only defined between
reduced marked graphs. 
\end{definition}


\section{Finding and factoring Whitehead moves}\label{sc:moves}

Recall that Whitehead moves suffice to connect any two reduced marked
graphs (roses) in Outer space \cite{ar:CV86}.  We generalize this result
to the setting of deformation spaces in Theorem \ref{th:whm} below. Then
Theorem \ref{th:moves} is proved by expressing Whitehead moves in terms
of slides, inductions, and $\As^{\pm 1}$--moves (Propositions
\ref{prop:whm-I} and \ref{prop:whm-II}). 

We need to introduce some notation. 
If $\Gamma,\Gamma'' \in \DS$ and $\Gamma =
\Gamma''_{F}$ for some forest $F \subset \Gamma''$ and $F_{0}$ is a
subforest of $F$, then we denote the marked graph $\Gamma''_{F -
F_{0}}$ by $\Gamma^{F_{0}}$. A forest $F_1 \subset \Gamma^{F_0}$ is
collapsible if and only if $(\Gamma^{F_0})_{F_1} \in \DS$, and when this 
occurs, we will abbreviate $(\Gamma^{F_0})_{F_1}$ as $\Gamma^{F_0}_{F_1}$. 
Given a graph of groups $\Gamma$ and an edge $e \in
E(\Gamma)$, the inclusion map $G_{e} \hookrightarrow G_{\iv(e)}$ may or
may not be surjective. We assign a \emph{label} to $e$, which is ``$=$'' if
$G_{e} \hookrightarrow G_{\iv(e)}$ is surjective and ``$\not=$'' otherwise. 

\begin{lemma}\label{lm:whm-ind}
Let $\Gamma,\Gamma' \in \DS$ be reduced marked graphs and
suppose that there is a marked graph $\Gamma'' \in \DS$ such that
$\Gamma = \Gamma''_{F}$ and $\Gamma' = \Gamma''_{F'}$ for non-trivial
finite forests $F,F' \subset \Gamma''$ that do not share an edge.  Then
there are edges $e \in E(F), e' \in E(F')$ such that one of the following
holds: 

\begin{enumerate}
    
    \item $\Gamma^{e}_{e'} \in \DS$ and $\Gamma^{e}_{e'}$ is
      reduced,  \label{c1}  
    
    \item there is an edge $f' \in E(\Gamma'')$ such that 
    $\Gamma^{e}_{e' \cup f'} \in \DS$ and $\Gamma^{e}_{e' \cup f'}$ is
    reduced,  \label{c2}  
        
    \item $\Gamma'^{e'}_{e} \in \DS$ and $\Gamma'^{e'}_{e}$ is reduced,
      or \label{c3}  
    
    \item there is an edge $f \in E(\Gamma'')$ such that 
    $\Gamma'^{e'}_{e \cup f} \in \DS$ and $\Gamma'^{e'}_{e \cup f}$ is
    reduced.  \label{c4}  
    
\end{enumerate}
\end{lemma}

Note that the lemma is symmetric in $\Gamma$ and $\Gamma'$. 

\begin{proof}
We begin with the following claim. 

\begin{claim}
If there are edges $e \in E(F)$, $e' \in E(F')$ such that $e'$ is collapsible
in $\Gamma^e$, then conclusion \eqref{c1} or \eqref{c2} holds. 
\end{claim}

\begin{proof}[Proof of Claim]
Replacing $e$ by $\bar{e}$ if necessary, we may assume that $e$ has label
$=$ in $\Gamma^e$. Then, since $\Gamma = \Gamma^{e}_{e}$ is reduced, every 
collapsible edge in $\Gamma^{e}$ must be incident to
$\iv(e)$. Furthermore, every such edge $f$ with $\iv(f) = \iv(e)$ has one
of three types: 
\begin{itemize}
    
    \item[] type 1: \  $\tv(f) \neq \tv(e)$ (which implies that 
    $\bar{f}$ has label $\not=$) 
    
    \item[] type 2: \  $\tv(f) = \tv(e)$ and 
    $\bar{f}$ has label $\not=$ 

    \item[] type 3: \  $\tv(f) = \tv(e)$ and 
    $\bar{f}$ has label $=$. 

\end{itemize}

Note that collapsing a type 2 edge always results in a reduced marked
graph. Also, after collapsing a type 1 edge, type 3 edges remain
collapsible and the other types become non-collapsible. Similarly, after
collapsing a type 3 edge, type 1 edges remain collapsible and the others
become non-collapsible. 

If $e'$ is of type 2 then conclusion \eqref{c1} holds. In fact, by the
observations above, the only way $\Gamma^e_{e'}$ can fail to be
reduced is if $\Gamma^e$ has collapsible edges $f_1$ of type 1
and $f_3$ of type 3, one of which is $e'$. Then $f_3$ is collapsible in
$\Gamma^e_{f_1}$, implying that $\Gamma^e_{f_1 \cup f_3} \in \DS$; and 
$\Gamma^e_{f_1 \cup f_3}$ is reduced, establishing \eqref{c2}. 
\end{proof} 

Returning to the lemma, we proceed by considering 
various configurations of the forests $F$ and $F'$.  Since $F \subseteq
\Gamma''$ is collapsible, in each component $F_{0}$ of $F$ there 
is a maximal subtree $F_{1} \subseteq F_{0}$ such that 
every (oriented) edge in $E(F_1)$ has label $=$ and every edge $e \in
E(F_0) - E(F_1)$ such that $\iv(e) \in V(F_0) - V(F_1)$, $\tv(e) \in
V(F_1)$ also has label $=$. We call $F_{1}$
the \emph{maximal stable subtree} of $F_{0}$.  

Let $e'$ be any edge of $F'$ with label $=$, and suppose that $e'$ does
not map to a loop in $\Gamma$.  Then there is a component $F_{0}$ of $F$
containing $\iv(e')$ but not $\tv(e')$.  Let $F_{1}$ be the maximal
stable subtree of $F_{0}$.  Notice that $\iv(e') \notin V(F_{1})$, since
otherwise $e'$ would be collapsible in $\Gamma$. 
Let $e$ be the first edge in the path in $F_0$ from $\iv(e')$ to
$F_1$. Note that collapsing $F - \{e, \bar{e}\}$ does not enlarge the vertex
group $G_{\iv(e')}$, and so $e'$ is collapsible in $\Gamma^e$. By the
Claim, conclusion \eqref{c1} or \eqref{c2} holds.  
Similarly, by symmetry, if there is an edge in $F$ with label $=$ which
does not map to a loop in $\Gamma'$, then conclusion \eqref{c3} or
\eqref{c4} holds. 

Therefore, we may assume that every edge of $F'$ maps to a loop
in $\Gamma$, and every edge of $F$ maps to a loop in $\Gamma'$. 
Now let $F_0$ be a component of $F$ and $F_1 \subseteq F_0$ its maximal
stable subtree. Choose a vertex $v \in V(F_1)$. There is an edge in $F_0$
incident to $v$, and there is a path in $F'$ joining the endpoints of
this edge (since it maps to a loop in $\Gamma'$). In particular, there is
an edge $e' \in E(F')$ with initial vertex $v$. Now let $\gamma$ be the path
in $F_0$ from $v$ to $\tv(e')$, which exists since $e'$ maps to a loop in
$\Gamma$. Let $e$ be the final edge of $\gamma$, with $\tv(e) =
\tv(e')$. Because $v \in V(F_1)$, collapsing $F - \{e, \bar{e}\}$ does
not enlarge the vertex groups at $v = \iv(e')$. 

The vertex group at $\tv(e')$ also remains unchanged. To see this, let
$F_2$ be the component of $F_0 - \{e, \bar{e}\}$ containing
$\tv(e)$. Since in $F_0$, $e$ separates $v$ from $F_2$, any vertex group
in $F_2$ is contained in the vertex group at $\tv(e) = \tv(e')$. 

Therefore $e'$ is collapsible in $\Gamma^e$, and by the Claim we are
finished. 
\end{proof}

\begin{theorem}\label{th:whm}
Any two reduced marked graphs in $\DS$ are related by a sequence of 
Whitehead moves.    
\end{theorem}

\begin{proof} 
Let $\Gamma, \Gamma' \in \DS$ be reduced marked graphs. First we consider
a special case, when $\Gamma$ and $\Gamma'$ satisfy the hypotheses of
Lemma \ref{lm:whm-ind}. In this case the theorem is proved by induction
on the number of edges in $F \cup F'$, as follows.  

Apply Lemma \ref{lm:whm-ind} and suppose that conclusion \eqref{c1}
holds. Then $\Gamma$ and $\Gamma^e_{e'}$ are related by a Whitehead
move. Also, $\Gamma''_{e'}$ collapses to $\Gamma^e_{e'}$ and to
$\Gamma'$, by collapsing the forests $F- \{e, \bar{e}\}$ and $F' - \{e',
\bar{e}'\}$ 
respectively. If both of these forests have no edges then $\Gamma^e_{e'}
= \Gamma'$ and we are done. Otherwise, since $\Gamma^e_{e'}$ and
$\Gamma'$ are reduced, both forests are non-trivial. Hence
$\Gamma^e_{e'}$ and $\Gamma'$ satisfy the hypotheses of Lemma
\ref{lm:whm-ind}, and by induction, they are related by Whitehead
moves. The cases \eqref{c2}, \eqref{c3}, \eqref{c4} are similar. 

For the general case, start with an elementary deformation from $\Gamma$
to $\Gamma'$, which can be written as 
\[\Gamma = \Gamma_0 \gets \Gamma_1 \to \Gamma_2 \gets \Gamma_3 \to \cdots
\gets \Gamma_{2n-1} \to \Gamma_{2n} = \Gamma'\] 
where each arrow is a sequence of collapses.
By inserting collape moves and their inverses at $\Gamma_{2i}$, we can
arrange that each $\Gamma_{2i}$ is reduced. 

Similarly, if the two forests being collapsed in $\Gamma_{2i} \gets
\Gamma_{2i+1} \to \Gamma_{2i+2}$ have shared edges, then the collapse
moves at these edges may be cancelled in pairs. Thus we may assume that
the two forests have no edges in common. Now $\Gamma_{2i}$ and
$\Gamma_{2i+2}$ satisfy the hypotheses of Lemma \ref{lm:whm-ind} for each
$i$, and are related by Whitehead moves, by the special case. 
\end{proof}


\begin{proposition}\label{prop:whm-I}
Any type I Whitehead move is a composition of slides and inductions, 
where the intermediate marked graphs are reduced. Moreover, the only
edge being slid over is one of the edges which is collapsed in the
Whitehead move. 
\end{proposition}

\begin{proof}
Let $\Gamma, \Gamma' \in \DS$ be reduced marked graphs such that
$\Gamma = \Gamma''_{e}$ and $\Gamma' = \Gamma''_{e'}$.  Orient $e,e'$
so that $\iv(e) = \iv(e')$.  Let $\{f_{\alpha}\}$ be the other
edges with initial vertex $\iv(e)$.
    
If $e \cup e'$ is not a cycle in $\Gamma''$ then $e$ and $e'$ both have
label $=$.  Then $\Gamma'$ is obtained from $\Gamma$ by sliding the
collection of edges $\{f_{\alpha}\}$ over $e'$.  There is no intermediate
marked graph in this case. 

Now suppose that $e \cup e'$ is a cycle and $e$ has label $=$ in
$\Gamma''$.  If $e'$ also has label $=$, then as above, $\Gamma'$ is
obtained from $\Gamma$ by sliding the collection of edges
$\{f_{\alpha}\}$ over the loop $e'$. Otherwise we can assume that
$\bar{e}, e'$ have label $\not=$ and $\bar{e}'$ has label $=$.  Then
$\Gamma'$ is obtained from $\Gamma$ by an induction move (in which the
same edge $e$ is expanded, but the edges $f_{\alpha}$ are left at
$\tv(e')$), followed by a slide of the collection of edges
$\{f_{\alpha}\}$ over the loop $\bar{e}$. Since $\Gamma'$ is reduced,
each $f_{\alpha}$ which is not a loop has label $\not=$ in $\Gamma''$.
This implies that the intermediate marked graph is reduced. 
\end{proof}

\begin{proposition}\label{prop:whm-II}
Any type II Whitehead move is a composition of slides and $\As$-- or
$\As^{-1}$--moves, where the intermediate marked graphs are reduced.
\end{proposition}

\begin{proof}
Let $\Gamma,\Gamma' \in \DS$ be reduced marked graphs such that 
$\Gamma = \Gamma''_{e}$ and $\Gamma' = \Gamma''_{e'\cup f'}$. This is the
situation of conclusion \eqref{c2} of Lemma \ref{lm:whm-ind}. The proof
of the Claim shows that $e'$ and $f'$ are of types 1 and 3. After
renaming, these edges of $\Gamma''$ must have the configuration shown in
Figure \ref{fig:lemma-fig2}. The labels are as shown because $\Gamma$ and
$\Gamma'$ are reduced. 
\begin{figure}[ht]
\psfrag{=}{{$=$}}
\psfrag{n}{{$\neq$}}
\psfrag{e}{{$e'$}}
\psfrag{f}{{$f'$}}
\psfrag{ep}{{$e$}}
\psfrag{g}{{$g$}}
\centering
\includegraphics[width=5cm]{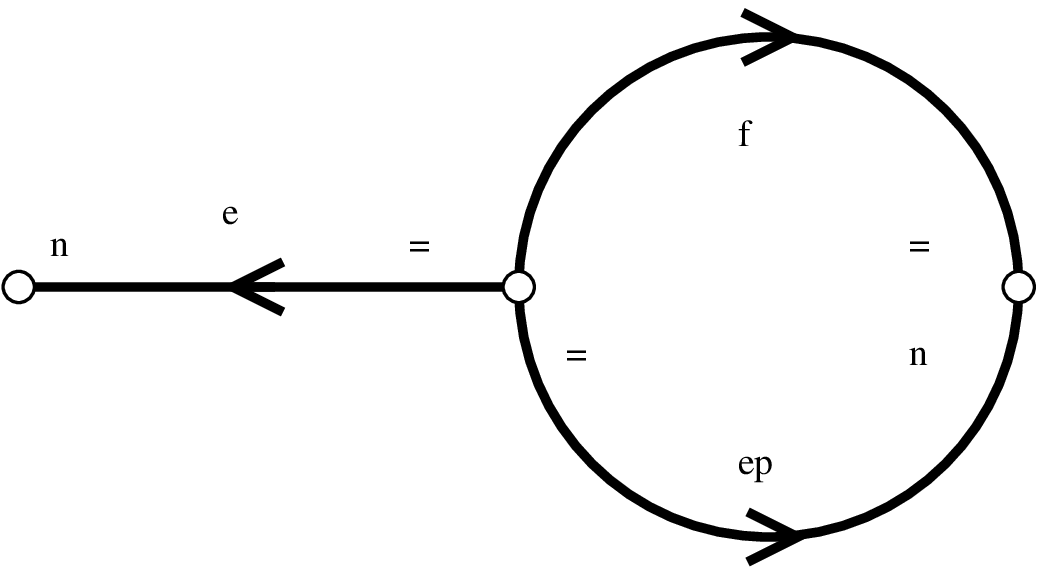}
\caption{ }
\label{fig:lemma-fig2}
\end{figure}
Let $\{g_{\alpha}\}$ be the other edges with initial vertex $\tv(e)$, and
$\{h_{\beta}\}$ the other edges with initial vertex $\iv(e)$. 
Now $\Gamma'$ is obtained from $\Gamma$ by first sliding the collection
$\{h_{\beta}\}$ over $e'$, then sliding $\{g_{\alpha}\}$ over $(\bar{f}',
e')$, and then performing an $\As^{-1}$--move. It is easy to verify 
that the intermediate marked graphs are reduced. 
\end{proof}

Theorem \ref{th:moves} now follows directly from Theorem \ref{th:whm} and
Propositions \ref{prop:whm-I} and \ref{prop:whm-II}. 

\medskip
The next result 
follows easily from Theorem
\ref{th:moves}, as explained in the introduction. The second statement is
included for use in \cite{ar:CF2}. 

\begin{corollary}\label{co:slides}
In a non-ascending deformation space of $G$--trees, any two reduced trees
are related by a finite sequence of slide moves, with all intermediate
trees reduced. Moreover, if $e$ is an edge of $T$ and a deformation from
$T$ to $T'$ never collapses $e$, then there is a sequence of slide moves
from $T$ to $T'$ in which no edge slides over $e$. 
\end{corollary}

\begin{proof}
Start with a deformation between reduced $G$--trees $T, T'$ in
$\DS$. There is a sequence of 
Whitehead moves joining $T$ to $T'$ by Theorem \ref{th:whm}, and these
moves are all of type I since $\DS$ is non-ascending (a type II
Whitehead move cannot occur by Proposition \ref{prop:whm-II}). Recall that 
in the proof of Theorem \ref{th:whm}, each time a type I Whitehead move
is factored out, the expansion and collapse comprising that move were
already present in the original elementary deformation. Thus, by
Proposition \ref{prop:whm-I}, there is a sequence of slide and induction
moves from $T$ to $T'$, and the only edges that are slid over were
expanded or collapsed in the original deformation. Lastly, there are no
induction moves since $\DS$ is non-ascending. 
\end{proof}

The rigidity theorem for $G$--trees \cite{ar:F02,ar:G03,ar:L05a} also
follows quickly from Theorem \ref{th:moves}. Recall that a $G$--tree $T
\in \DS$ is \emph{rigid} if it is the only reduced $G$--tree in $\DS$.

\begin{corollary}[Levitt]\label{co:rigid}
A $G$--tree that is not the Bass--Serre tree of an ascending
HNN-extension is rigid if and only if, for any two edges $e,f$ such 
that $\iv(e) = \iv(f) = v$ and $G_{e} \subseteq G_{f}$, one of the 
following conditions holds:
\begin{enumerate}
    \item $e \in G f$, \label{con1} 
    
    \item $e \in G \bar{f}$ and $G_{e} = G_{f}$, or \label{con2} 
    
    \item there is an edge $f'$ such that $\iv(f') = v$, $f' 
    \in G \bar{f}$ and $G_{f} = G_{f'} = G_{v}$ and there are only 
    three $G_{v}$--orbits of edges at $v$. \label{con3} 
\end{enumerate} 
\end{corollary}

\begin{proof}
By Theorem \ref{th:moves} it is clear that such a $G$--tree is rigid
if and only if it does not admit a slide, induction or $\As^{\pm
1}$--move resulting in a different $G$--tree.  Given $e,f$ as above, if
$e \notin Gf \cup G\bar{f}$ then there is a (possibly trivial) slide move
of $e$ over $f$.  This slide move is trivial only under
the conditions of \eqref{con3}. 
If $e \in G\bar{f}$ and $G_{e} \neq G_{f}$ then the $G$--tree admits an
$\As^{\pm 1}$--move or induction move. (The image of $e$ in the quotient
marked graph is a loop; either it is an ascending loop, or it satisfies
the criteria for an $\As$--move given in Definition \ref{def:A}.) 
Note that an induction move is
always non-trivial, except possibly in the case of an ascending HNN
extension.  Thus, if $T$ is rigid, then $e$ and $f$ satisfy one
of the three conditions. 
For the converse, if $T$ admits a slide, induction
or $\As^{\pm 1}$--move, then there is a pair of edges $e,f$ that do not
satisfy any of the three conditions.
\end{proof}

Theorem \ref{th:moves} implies that a $G$--tree that is the
Bass--Serre tree of an ascending HNN-extension is rigid if and only if
it does not admit a nontrivial induction move.  See
\cite[Theorem~2]{ar:L05a} for algebraic conditions on the monodromy
$\phi\co A \hookrightarrow A$ characterizing when the $G$--tree does not
admit a nontrivial induction move.


\def\cprime{$'$}

\end{document}